\documentclass[12pt]{amsart}
\usepackage[all]{xy}

\swapnumbers
\newtheorem{thm}[subsection]{Theorem}
\newtheorem{prop}[subsection]{Proposition}
\newtheorem{lemma}[subsection]{Lemma}
\newtheorem{cor}[subsection]{Corollary}
\newtheorem{rem}[subsection]{Remark}

\def\oo{\omega}
\def\PP{{\mathcal P}}

\def\arbreA{\xymatrix@R=3pt@C=3pt{
&& \\
&*{}\ar@{-}[ur] \ar@{-}[ul] \ar@{-}[d]     &\\
&&
}}

\def\arbreBA{\xymatrix@R=2pt@C=2pt{
&&&&\\
&&&*{}\ar@{-}[ul] & \\
&&*{}\ar@{-}[uurr] \ar@{-}[uull] \ar@{-}[d]     &&\\
&&&&
}}

\def\arbreAB{\xymatrix@R=2pt@C=2pt{
&&&&\\
&*{}\ar@{-}[ur] &&& \\
&&*{}\ar@{-}[uurr] \ar@{-}[uull] \ar@{-}[d]     &&\\
&&&&
}}

\def\arbreBB{\xymatrix@R=2pt@C=2pt{
&&&&\\
&&&& \\
&&*{}\ar@{-}[uurr] \ar@{-}[uull] \ar@{-}[d] \ar@{-}[uu]     &&\\
&&&&
}}

\def\arbreABC{\xymatrix@R=1pt@C=1pt{
&&&&&&\\
&*{}\ar@{-}[ur] &&&&& \\
&&*{}\ar@{-}[uurr] &&&&\\
&&&*{}\ar@{-}[uuurrr] \ar@{-}[uuulll] \ar@{-}[d] &&&\\
&&&&&&
}}

\def\arbreBAC{\xymatrix@R=1pt@C=1pt{
&&&&&&\\
&&&*{}\ar@{-}[ul] &&& \\
&&*{}\ar@{-}[uurr] &&&&\\
&&&*{}\ar@{-}[uuurrr] \ar@{-}[uuulll] \ar@{-}[d] &&&\\
&&&&&&
}}

\def\arbreACA{\xymatrix@R=1pt@C=1pt{
&&&&&&\\
&*{}\ar@{-}[ur] &&&&*{}\ar@{-}[ul] & \\
&&&&&&\\
&&&*{}\ar@{-}[uuurrr] \ar@{-}[uuulll] \ar@{-}[d] &&&\\
&&&&&&
}}

\def\arbreCAB{\xymatrix@R=1pt@C=1pt{
&&&&&&\\
&&&*{}\ar@{-}[ur] &&& \\
&&&&*{}\ar@{-}[uull] &&\\
&&&*{}\ar@{-}[uuurrr] \ar@{-}[uuulll] \ar@{-}[d] &&&\\
&&&&&&
}}

\def\arbreCBA{\xymatrix@R=1pt@C=1pt{
&&&&&&\\
&&&&&*{}\ar@{-}[ul] & \\
&&&&*{}\ar@{-}[uull] &&\\
&&&*{}\ar@{-}[uuurrr] \ar@{-}[uuulll] \ar@{-}[d] &&&\\
&&&&&&
}}

\def\arbreBAdecore{\xymatrix@R=2pt@C=2pt{
&&&&\\
&&&*{}\ar@{-}[ul] &\epsilon(i+1) \\
&&*{}\ar@{-}[uurr] \ar@{-}[uull] &&\\
&&\epsilon(i)&&
}}

\def\arbreABdecore{\xymatrix@R=2pt@C=2pt{
&&&&\\
\epsilon(i)&*{}\ar@{-}[ur] &&& \\
&&*{}\ar@{-}[uurr] \ar@{-}[uull]      &&\\
&&\epsilon(i+1)&&
}}

\begin{document}

\author[J.-L. Loday]{Jean-Louis Loday}
\address{Institut de Recherche Math\'ematique Avanc\'ee\\
    CNRS et Universit\'e Louis Pasteur\\
    7 rue R. Descartes\\
    67084 Strasbourg Cedex, France}
\email{loday@math.u-strasbg.fr}
\urladdr{www-irma.u-strasbg.fr/{$\sim$}loday/}
\title{Inversion of integral series enumerating  planar trees}
\subjclass[2000]{}
\keywords{Integer sequence, generating series, planar tree, excluded pattern, operad, Koszul duality}

\date{\today}

\begin{abstract}
We consider an integral series $f(X,t)$ which depends on the choice of a set $X$ of labelled planar rooted trees. We
prove that its inverse for composition is of the form $f(Z,t)$ for another set $Z$ of trees, deduced from
$X$. The proof is self-contained, though inspired by the Koszul duality theory of quadratic operads.
\end{abstract}

\maketitle

\section{Introduction} \label{S:int}

Let $I$ be a finite set of indices. Let $Y_n\times I^n$ be the set of planar binary rooted trees whose $n$ vertices are
 labelled by elements in the index set $I$. Let $X$ be a subset of $Y_2\times I^2$ and let $Z$ be its complement.
Define $X_n$ as the subset of $Y_n\times I^n$ made of labelled trees whose local patterns are in $X$. In other
words, a tree is in $X_n$ if for every pair of adjacent vertices the subtree defined by this pair is in $X$. By
convention
$X_0=Y_0\times I^0$ and $X_1=Y_1\times I$. From the definition of $X_n$ it comes immediately $X_2=X$.
The set $Z$ determines similarly a sequence $Z_n$.

The alternate generating series of $X$ is by definition
\begin{equation*}f(X,t) := \sum_{n\geq 0} (-1)^{n+1} (\# X_{n}) t^{n+1} = -t + (\# I)t^2 -  (\# X)t^3 +\cdots \ .
\end{equation*}
\noindent {\bf Theorem.} {\it If $Z$ is the complement of $X$, i.e.~$X \sqcup Z= Y_2\times I^2$, then the generating series of $X$ and $Z$ are inverse to each other for composition:
$$f(X, f(Z,t)) = t .$$
}

For some choices of $I$ and $X$ the integer sequence $(\# X_{n})_{n\geq 0}$ appear in the data base ``On-line Encyclopedia of Integer sequences!" \cite {Sloane},
but for some others they do not.

 Here is an application of this theorem.
Given an integer sequence ${\underline a}=(a_0, \ldots, a_n, \ldots)$ it is often interesting to know a combinatorial
 interpretation of these numbers, that is to know a family $X_n$ of combinatorial objects such that $a_n = \# X_n$.
The theorem provides a solution for some integer sequences as follows. Suppose that the inverse for composition
of the alternate series of ${\underline a}$ gives an integer sequence ${\underline b}$ which can be interpreted
combinatorially by labelled trees. Then the integer sequence ${\underline a}$ admits also such an interpretation.

Our proof of the theorem consists in constructing a chain complex whose Poincar\'e series is exactly
 $f(X, f(Z,t))$. Then we prove that this chain complex is acyclic (i.e.~the homology groups are 0 except $H_1$ which
is of dimension 1) by reducing it to the sum of subcomplexes which
turn out to be augmented chain complexes of standard simplices. Hence the Poincar\'e series is $t$. 

Our proof is self-contained but the idea of considering this particular chain complex is inspired by the theory of quadratic operads. Indeed the
choice of
$X$ determines a certain type of algebras, i.e.~a certain quadratic operad, and the choice of $Z$ gives the ``dual
operad" in the Koszul duality sense cf.~\cite{G-K}. Then the chain complex is the Koszul complex attached to this dual  pair of
operads. So our main theorem gives a large family of Koszul operads.

We give all the details for the case of binary trees, but this method can be generalized to planar trees. We outline the  case of $k$-ary trees in the last section. A surprizing consequence is the following property of the Catalan numbers $c_n$. The series $h(t) = \sum_{n\geq 0}(-1)^{n+1}c_nt^{3n+1}$ is its own inverse for composition: $h(h(t))=t$.

After the release of the first version of this paper, I was informed by Prof.~I.~Gessel that his student S.F.~Parker obtained the same result by combinatorial methods in her thesis (unpublished). A far reaching generalization of our result has been obtained subsequently by R.~Bacher in \cite{B} using different techniques.

\section{Labelled trees}
\subsection{Planar binary rooted trees}
Denote by $Y_n$ the set of {\it planar binary rooted trees} of degree $n$, that is with $n$ vertices (with valence at least $2$):
\begin{displaymath}
Y_0 = \{\  \vert \ \}\  , \quad Y_1=\big\{ \vcenter{\arbreA}\big\}\  , \quad 
{Y_2=\Big\{ L:=\vcenter{\arbreBA}}\ ,\  {R:= \vcenter{\arbreAB }}\Big\}
\end{displaymath}
\begin{displaymath}
{Y_3=\Bigg\{ \vcenter{\arbreABC}}\ ,\  {{}\vcenter{\arbreBAC }},\  {{}\vcenter{\arbreACA }},\  {{}\vcenter{\arbreCAB }},\  {{}\vcenter{\arbreCBA }}\Bigg\}
\end{displaymath}

The number of elements in $Y_n$ is the so-called {\it Catalan number} $c_n=\frac {(2n)!}{n!(n+1)!} $ cf. \ref{Examples} (b). Let $I$ be a
 finite set of indices. By definition a {\it labelled tree} is a planar binary rooted tree such that
each internal vertex is labelled by an element of $I$. These elements need not be distinct. Therefore the set of
labelled trees of degree $n$ is in bijection with $Y_n\times I^n$.

An element of $Y_2\times I^2$ is either of the form $(L; i_1,
i_2)$ or of the form $(R; i_1, i_2)$.

Let $X$ be a subset of $Y_2\times I^2$ and let $Z$ be its complement. We define a subset $X_n$ of  $Y_n\times I^n$ as
 follows.  A pair of adjacent vertices in the labelled tree $y$ determines a subtree of degree $2$, called a
\emph {local pattern}. The labelled tree $y$ is in $X_n$ if and only if all its local patterns belong to $X$. In
other words we exclude all the trees which have a local pattern which belongs to $Z$. It is clear that $X_2=X$. By
convention we define
$X_0=Y_0\times I^0$ and $X_1=Y_1\times I$. 

The alternate generating series of $X$ is determined by the integer sequence $(\# X_n)_{n\geq 0}$ as follows:
\begin{equation*}f(X,t) := \sum_{n\geq 0} (-1)^{n+1} (\# X_{n}) t^{n+1} = -t + (\# I)t^2 -  (\# X)t^3 +\ldots \ .
\end{equation*}

If, instead of $X$, we start we $Z$, then we get another family  $Z_n$. For $n=2$, $Z_2$ is the complement of $X_2$, but this property does not hold for higher $n$'s.

\begin{thm} \label{THM}  If $Z$ is the complement of $X$, i.e.~$X \sqcup Z= Y_2\times I^2$, then the generating series
of
$X$ and
$Z$ are inverse to each other for composition:
$$f(X, f(Z,t)) = t .$$
\end{thm}
The proof is given in the next section.

\subsection{Examples}\label{Examples}  We list a few interesting examples of integer sequences and their dual which appear in the study of quadratic operads (cf. ~\cite{L2}). The notation is as follows: the sequence $(a_0, \cdots , a_n, \cdots )$ is such that $a_n = \#X_n$ and $f(t)= \sum_{n\geq 0}(-1)^{n+1}a_nt^{n+1}$.
The dual sequence is  $(b_0, \cdots , b_n, \cdots )$ where $b_n = \#Z_n$ and $g(t)= \sum_{n\geq 0}(-1)^{n+1}b_nt^{n+1}$.

Observe that $a_0 = 1$, $a_1 = \# I$, $a_2 = \#X$ and $0\leq a_2\leq 2(a_1)^2$. In the following examples  we can write $g(t)$ as a rational function, hence we get a
 combinatorial interpretation of the integer sequence $\underline a$ whose alternate series $f(t)$ is determined by the functional equation $g(f(t))=t$.
\\

\noindent (a) $(1, 1,1, \ldots, 1, \ldots )$ versus itself.
\begin{itemize}
\item $I=\{1\}$, $X= \{L\}$ and $Z=\{R\}$. 
\item $f(t)=g(t) = \frac {-t}{1+t}$.
\end{itemize}

\noindent (b) $(1, 1, 2, 5, 14, 42, 132, \cdots, c_n, \ldots )$ versus $(1,1,0, \ldots , 0, \ldots $).
\begin{itemize}
\item  $I=\{1\}, Z=\emptyset$,  $X=Y_2$.
\item $g(t) = -t+t^2$.
\item We get the well known functional equation for the generating series of the Catalan numbers $c(t):= \sum_{n\geq 0} c_nt^n$,
\[ tc(t)^2 - c(t) +1 = 0. \]
\end{itemize}

\noindent (c)  $(1, 2, 6, 22, 90, \cdots, 2C_{n}, \ldots )$ versus $(1,2,2, \ldots , 2, \ldots $).
\begin{itemize}
\item$C_n$ is the super Catalan number (also called Schr\"oder number), that is the number of planar trees with $n+1$ leaves.
\item $I=\{1,2\}$, $Z=\{(L;1,1), (R;2,2)\}$ and $X$ has 6 elements. 
\item It is immediate to see that $Z_n$
has only two elements: the right comb indexed by $1$'s and the left comb indexed by $2$'s. Therefore
$g(t) = \frac {-t+t^2}{1+t}$. On the other hand one can show that, for $n\geq
2$, there is a bijection between the elements of $X_n$ and two copies of the set of planar trees with $n+1$ leaves (see \cite{L-R2} for a variant of this result). The theorem gives the well known functional equation for
the generating series of the super Catalan numbers $C(t):= \sum_{n\geq 0} C_nt^n$,
\[ tC(t)^2 +(1-t)C(t) -1 = 0. \]
\end{itemize}

\noindent (d)  $(1, 2, 6, 21,80, \ldots )$ versus $(1,2,2,1,0,\ldots, 0 , \ldots )$.
\begin{itemize}
\item $I=\{1,2\}$, $Z=\{(L,2,1), (R,1,2)\}$ and $X$ has 6 elements. 
\item It is immediate to see that $X_3$
has only one element and that $X_n$ is empty for $n\geq 4$. Hence $g(t)=-t+2t^2-2t^{3}+t^4$.
\item This example and the previous one show that the integer sequence determined by $X$ does not depend only on the number of elements of $I$ and $X$.
\end{itemize}

\noindent (e)  $(1, 2, 7, 31, 154, \ldots )$ versus $(1,2,1,1,\ldots, 1 , \ldots )$.
\begin{itemize}
\item Let $I=\{1,2\}$, $Z=\{(L;1,1)\}$ and $X$ has 7 elements.  
\item It is immediate to see that for $n\geq 2$, 
$Z_n$ has only one element: the right comb indexed by $1$'s. Hence  $g(t)= \frac {-t+t^2+t^3}{1+t}$. 
\end{itemize}

\noindent (f)  $(1, 3, 17, 121, 965, \cdots )$ versus $(1,3,1,1,\ldots, 1 , \ldots )$, and \\
\indent $(1, k, 2k^2-1, 5\, 3^k -5k+1\cdots, ?, \cdots )$ versus $(1,k,1,1,\ldots, 1 , \ldots )$.
\begin{itemize}
\item $I=\{1,\cdots , k\}$, $Z=\{(L;1,1)\}$ and $X$  has $2k^2-1$ elements.
\item It is immediate to see that, for $n\geq 2$,  $Z_n$ has only one element.
Hence $g(t)= \frac {-t+(k-1)(1+t)t^2}{1+t}$. 
\end{itemize}


\noindent (g)  $(1, 3, 14, 80, 510, \cdots )$ versus $(1,3,4,5,\ldots, n+2,\ldots )$.
\begin{itemize}
\item I=$\{1,2,3\}$, $Z= \{(L;1,1), (L;2,1), (R;2,2), (R;3,3)\}$, $X$ has $14$ elements.
\item One checks that $Z_n$ is made of $(n+1)+1$ elements. Hence $g(t)=\frac{t(-1+t+t^2)}{(1+t)^2}$.
\end{itemize}

\noindent (h)  $(1, 4, 23, 156, 1162, \cdots )$ versus $(1,4,9,16,\ldots, (n+1)^2,\ldots )$.
\begin{itemize}
\item I=$\{\nwarrow, \nearrow, \searrow, \swarrow \}$, in the following we write $Rij$ in place of $(R;i,j)$:
\[ Z= \left\{
\begin{array}{ccc}
R\nwarrow \nwarrow  & R\nearrow \nwarrow  & L\nearrow \nearrow   \\
R\swarrow \nwarrow  & R\searrow \nwarrow  & L\searrow \nearrow  \\
R\swarrow \swarrow  & L\searrow \swarrow  & L\searrow \searrow   \\
\end{array}\right\}
\]
\item We have shown in \cite{A-L} Proposition 4.4 that  $Z_n= (n+1)^2$, hence $g(t)=\frac{t(-1+t)}{(1+t)^3}$.
\item The first sequence has been given a different combinatorial interpretation in terms of connected non-crossing configurations in \cite{Flajolet}.
\end{itemize}

\noindent (i)  $(1, 9, 113, \cdots )$ versus $(1,9,49,\ldots,  )$.
\begin{itemize}
\item $I$ is a set of $9$ indices denoted
$ \begin{array}{ccc}
\nwarrow  & \uparrow  & \nearrow  \\
 \leftarrow  & \circ  & \rightarrow  \\
 \swarrow  & \downarrow  & \searrow  \\
\end{array}
$
and $Z$ is made of the following 49 elements (indexed by the cells of $\Delta^2\times \Delta^2$):

\[ \begin{array}{ccccccc}
R\nwarrow \nwarrow  & R\nearrow \nwarrow  & L\nearrow \nearrow  & R\nearrow \uparrow  & R\nwarrow \uparrow  & R \uparrow \nwarrow & R\uparrow \uparrow \\
R\swarrow \nwarrow  & R\searrow \nwarrow  & L\searrow \nearrow  & R\searrow \uparrow  & R\swarrow \uparrow  & L\downarrow \nwarrow  & R\downarrow \uparrow  \\
R\swarrow \swarrow  & L\searrow \swarrow  & L\searrow \searrow  & R\searrow \downarrow  & L\swarrow \downarrow  & L\downarrow \swarrow  & L\downarrow \downarrow  \\
R\swarrow \leftarrow  & R\searrow \leftarrow  & L\searrow \rightarrow  & R\searrow \circ  & R\swarrow \circ  & L\downarrow \leftarrow  & R\downarrow \circ  \\
R\nwarrow \leftarrow  & R\nearrow \leftarrow  & L\nearrow \rightarrow  & R\nearrow \circ  & R\nwarrow \circ  & L\uparrow \leftarrow  & R\uparrow \circ  \\
R\leftarrow \nwarrow  & L\rightarrow \nwarrow  & L\rightarrow \nearrow  & R\rightarrow \uparrow  & L\leftarrow \uparrow  & L\circ \nwarrow  & L\circ \uparrow  \\
R\leftarrow \leftarrow  & R\rightarrow \leftarrow  & L\rightarrow \rightarrow  & R\rightarrow \circ  & R\leftarrow \circ  & L\circ \leftarrow  & R\circ \circ  \\
\end{array}
\]

\item Unfortunately we do not know how to compute the number of elements in $Z_n$. This example is strongly related to dendriform trialgebras \cite{L-R1} and motivated by the  ennea-algebras \cite {Leroux}.
\end{itemize}

Added after the release of the first version: a thorough study of this example has been performed in \cite{B}.

\section{Koszul complex and the Theorem}

\subsection{Koszul complex}

Given a planar binary rooted tree $y\in Y_n$ one  numbers the leaves from left to
right by $0, \ldots , n$. Accordingly, one numbers the vertices by $1, \ldots, n$, the $i$th
vertex being in between the leaves $i-1$ and $i$. So a decoration is a map $\epsilon$ from
$\{1,\ldots, n\}$ to $I$. A vertex of $y$ is said to be a \emph{cup} if it is directly
connected to two leaves (no intermediate vertex). In the following example $z$ has a cup at
vertex $1$ and vertex $3$:
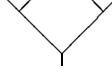
\[ \arbreACA \]
The grafting of two treees $y$ and $y'$ is the new tree $y\vee y'$ obtained from $y$ and $y'$ by joining the roots to a new vertex and adding a new root. For instance the above tree is the grafting $\vcenter{\arbreA} \vee \vcenter{\arbreA}$.

We define a chain complex ${\mathcal K}_* = ({\mathcal K}_n , d)_{n>0}$ over the field ${\mathbb K}$ as follows. The
space of $(n+1)$-chains is
$$ {\mathcal K}_{n+1} := \bigoplus {\mathbb K}[Z_n \times X_{i_0} \times \cdots \times  X_{i_n}]$$
where the sum is extended to all $(n+1)$-tuples $(i_0, \cdots , i_n)$, where $i_j\geq 0$. The boundary map $d :
{\mathcal K}_{n+1}\to {\mathcal K}_{n}$ is of the form $d= \sum_{i=1}^n (-1)^i d_i$, where
$d_i$ sends a basis vector to a basis vector or $0$ according to the following rule.

Let $z\in Z_n$ and $x_j\in X_{i_j}$. If the $i$th vertex of $z$ is not a cup, then
$d_i(z;x_0,\ldots , x_n):=0$. If the $i$th vertex of $z$ is a cup, then 
$$d_i(z;x_0,\ldots
, x_n):= (d_i(z);x_0,\ldots , x_{i-1}\vee_{\epsilon (i)} x_{i}, \ldots , x_n)$$
where $d_i(z)$ is the labelled tree obtained from $z$ by deleting the $i$th vertex
(replace it by a leaf), and where $\vee_{\epsilon (i)} $ means the grafting with ${\epsilon
(i)}$ as the decoration of the new vertex. If it happens that the labelled tree $x_{i-1}\vee_{\epsilon (i)} x_{i}$ contains a pattern in $Z$, then we put $d_i(z;x_0,\ldots , x_n):=0$.

\begin{lemma} $d^2=0$.
\end{lemma}
\noindent Proof. Let $\omega=(z;x_0,\ldots, x_n)$. It is sufficient to prove that $d_id_j =
d_{j-1}d_i$ for $i<j$ (presimplicial relation). If $i<j+1$, then the actions of $d_i$ and
$d_j$ on $\oo$  are sufficiently far apart so that they commute (the indexing $j-1$
comes from the renumbering). In the case $j=i+1$ we will prove that $d_id_{i+1} (\omega)=0=
d_{i}d_i (\omega)$. We are, locally in $z$, in one of the following two situations:
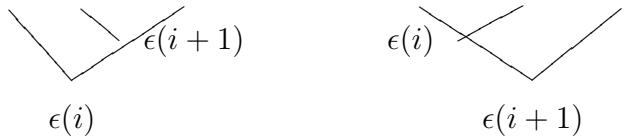
\[ \arbreBAdecore \qquad \qquad \arbreABdecore\]
In the first situation $d_id_i (\omega)=0$ because the $i$th vertex is not a cup. If $i+1$ is not a cup, then 
$d_id_{i+1} (\omega)=0$ because $d_{i+1} (\omega)=0$. If $i+1$ is a cup, then 
$d_id_{i+1} (\omega)=0$ because one of the entries of $d_id_{i+1} (\omega)$ is $a\vee_{\epsilon (i)}
(b\vee_{\epsilon (i+1)} c)$ which has a local pattern in  $Z$ and so is $0$ in $X_l$. 

The proof
is similar in the second situation.\hfill $\square$
\\

The chain complex ${\mathcal K}_*$ is called the \emph{Koszul complex} of $X$ (see
section~\ref{S:op} for an explanation of this terminology).

\subsection{Extremal elements}
By definition a basis vector $\omega=(z;x_0,\ldots, x_n)$ of  ${\mathcal K}_n$ is an \emph{extremal element} 
if there does not exist a basis vector $\omega '$ such that $d_i(\omega ') = \omega$ for
some $i$. 

\begin{prop}\label{extremal}
For each extremal element $\omega$ with $k$ cups, the basis vectors $d_{i_1}\cdots 
d_{i_r}\omega$ span a subcomplex ${\mathcal
K}_{\omega}$  of $\mathcal K$ which is isomorphic to the augmented chain complex of the
standard simplex $\Delta ^{k-1}$.
\end{prop}

\noindent Proof.  Let $\omega=(z;x_0,\ldots, x_n)$ be an extremal element.  The graded subvector space of ${\mathcal K}_*$ spanned by the
elements $d_{i_1}\cdots 
d_{i_r}\omega$ is stable by $d$ and so forms a subcomplex.

Let us now prove the isomorphism. We claim that $d_{i_1}\cdots  d_{i_r}\omega$ is
non-zero if and only if the indices
$i_j$ are such that the vertices ${i_j}$ are cups. Indeed the ``only if"
case is immediate. In the other direction: if $d_{i_1}\cdots  d_{i_r}\omega\ne 0$, then this would say that there is an $l$, such that the vertex $l$ is  a cup and 
$d_l(\oo)=(d_l(z); \ldots, (a\vee_u b)\vee_v c, \ldots)$ with $(R; u,v)$   in $Z$, or $d_l(\oo)=(d_l(z); \ldots, a\vee_u( b\vee_v c, \ldots)$ with $(L; u,v)$   in $Z$. So could construct ${\widehat \oo}= (\widehat z; \ldots, a, b, c, \ldots)$ so that 
$d_l({\widehat \oo})= \oo$ and $\omega$ would not be extremal.

We construct a bijection between the cells of $\Delta ^{k-1}$ and the set of non-zero vectors $\{
d_{i_1}\cdots  d_{i_r}\omega\}$ by sending the $j$th vertex  (number $j-1$) of $\Delta^{k-1}$ to  
$d_{i_1}\ldots \widehat{d_{i_j}}\ldots d_{i_k}\omega$ where $i_j$ is the $j$th cup of $x$. It is
immediate to verify that the boundary map in the chain complex of the standard simplex
corresponds to the boundary map $d$ by this bijection. Observe that, under this bijection, the
big cell of the simplex is mapped to $\oo$ and the generator of the augmentation space is mapped to  $d_{i_1}\ldots   d_{i_k}\omega$.\hfill $\square$

\begin{prop}\label{sum}
The chain complex ${\mathcal K}_*$ is isomorphic to $\bigoplus _{\omega}{\mathcal K}_
{\omega}$ where the sum is taken over all the extremal elements $\omega$.
\end{prop}

\noindent Proof. Let us show that any basis vector belongs to ${\mathcal K}_
{\omega}$ for some extremal element $\omega$. If $\omega$ is extremal, then the proposition holds. If
not, then there exists an element ${\omega}_1$ such that $d_i({\omega}_1)=\omega$ for
some
$i$, and so on. The process stops after a finite number of steps because $(z;\vert , \ldots, \vert )$ is extremal.

Now it is sufficient to prove that, if a basis vector belongs to ${\mathcal K}_{\omega}$ and 
to ${\mathcal K}_{\omega'}$, then $\omega = \omega '$.

Let $\omega=(z;x_0,\ldots, x_n)$ and $\omega'=(z';x_0',\ldots, x_n')$. If
$d_i(\oo)=d_j(\oo')\ne 0$, then $i$ is a cup of $z$, $j$ is a cup of $z'$ and $d_i(z)=d_j(z')$. If
$i< j$, then there exists $\bar\oo$ such that $d_j(\bar \oo) = \oo$, $d_i(\bar \oo)=\omega '$, and
therefore $\oo$ is not extremal. So we have $i=j$. 

If $d_i(\oo)=d_i(\oo')\ne 0$, then it is of the form $(\bar z; \ldots, a\vee_{\epsilon (i)} b,
\ldots )$. But the element $(z; \ldots, a, b,\ldots )$ where $z$ is the labelled tree obtained from
$\bar z$ by replacing the $i$th leaf by a cup and putting $\epsilon(i)$ as a decoration, is
the only element such that $d_i(\oo)=(\bar z; \ldots, a\vee_{\epsilon (i)} b,
\ldots )$. Hence $\oo = \oo '$.

So we have proved that any basis vector belongs to one and only one subcomplex of the form ${\mathcal K}_
{\omega}$. \hfill$\square$

\begin{rem} \end{rem}In order to visualize these two proofs it is helpful to think of the element $\omega=(z;x_0,\ldots, x_n)$ as a single graph (with a ``horizon") obtained by gluing the $x_j$'s to the leaves of $z$. The horizon indicates where to cut to get the $x_j$'s back. The operator $d_j$, where $j$ is the number of a cup, consists in lowering the horizon under the relevant vertex.

\begin{cor}\label{acyclicity}
For any choice of $X$ the Koszul complex ${\mathcal K}_*$ is acyclic.
\end{cor}
\noindent Proof. By Propositions~\ref{sum} and \ref{extremal} the homology of the Koszul complex is trivial since the standard simplex
is contractible. There is only one exception in dimension $1$ since the subcomplex
corresponding to the extremal element $\oo = (\vert ;\vert )$ is ${\mathbb K}$  in dimension $1$. So we have $H_n({\mathcal
K}_*) = 0$ for $n>1$ and $H_1({\mathcal
K}_*) = {\mathbb K}$. \hfill $\square$

\begin{prop}\label{Pseries}
The  Poincar\'e series of the Koszul complex ${\mathcal K}_*$ is equal to $f(Z,f(X,t))$.
\end{prop}

\noindent Proof. Let us call $w= n + i_0 + \cdots + i_n$ the \emph{weight} of an element $\oo\in Z_n\times X_{i_0}\times
\cdots \times X_{i_n}$. From the definition of $d_i(\oo)$ we see that the weight of $d_i(\oo)$ is also $w$. Therefore
the Koszul complex is the direct sum of subcomplexes ${\mathcal K}_*^{(w)}$ made of all the elements of weight
$w$. For a fixed weight $w$ the complex ${\mathcal K}_*^{(w)}$ is finite, beginning with ${\mathbb K}[Z_w\times (X_0)^{w+1}]$,
ending with
${\mathbb K}[Z_0\times X_w]$. More generally one has ${\mathcal K}_n^{(w)}=\bigoplus {\mathbb K}[Z_n\times X_{i_0}\times \cdots
\times X_{i_n}]$ where the sum is extended over all the $(n+1)$-tuples $(i_0,\ldots , i_n)$ such that $n + i_0 +
\cdots +i_n = w$.

Let $a_n:= \#X_n$ and $b_n := \#Z_n$ so that $f(X,t)=  \sum_{n\geq 1} (-1)^{n+1} a_n t^{n+1}$ and  
 $f(Z,t)= \sum_{n\geq 1} (-1)^{n+1} b_n t^{n+1}$. From the explicit description of
${\mathcal K}_n^{(w)}$ we check that the Euler-Poincar\'e characteristic of
${\mathcal K}_*^{(w)}$ is precisely the coefficient of $(-1)^wt^{w+1}$ in the expansion of 
$$ \sum_{n\geq 1} (-1)^{n+1} b_n \big(  \sum_{m\geq 1} (-1)^{m+1} a_m t^{m+1}\big)^n\ .$$
Therefore the Poincar\'e series $\sum_{w\geq 0}(-1)^{w}\chi({\mathcal K}_*^{(w)})t^{w+1}$ is equal to\\
$f(Z,f(X,t))$.
\hfill
$\square$

\subsection{End of the proof of Theorem~\ref{THM}}
By Proposition~\ref{Pseries} it suffices to show that the Poincar\'e series of ${\mathcal K}_*$ is $t$.
The Poincar\'e series of a complex is the same as the Poincar\'e series of its homology. Since the homology of 
${\mathcal K}_*$ is 0, except in weight 0 where it is ${\mathbb K}$ by Corollary \ref{acyclicity}, the Poincar\'e series is $t$. \hfill $\square$

\section{Operadic interpretation.}\label{S:op}
\subsection{Algebraic operad} Let $I$ be a finite set of indices,  $X$ be a subset of $Y_2\times I^2$ and $Z$  its complement. 
Over the field ${\mathbb K}$ we define a type of algebras, denoted $\PP$, as follows. There is one
binary operation $\circ_i$ for any $i\in I$ and the relations are 
$$(x\circ_i y)\circ_j z = 0 \hbox{ if } (R;i,j)\in Z\hbox{ and } 
x\circ_i ( y\circ_j z )= 0 \hbox{ if } (L;i,j)\in Z.
$$
It is immediate to check that the free algebra of type $\PP$ on one generator admits
$X_{n-1}$ as a basis of the homogeneous part of degree $n$, $n\geq 1$. The generator is
the unique element of $X_0$, that is $\vert $ . So the operad $\PP$ determined by this type
of algebras is such that $\PP(n)={\mathbb K}[X_{n-1}]\otimes {\mathbb K}[S_n]$, where $S_n$ is the symmetric
group. In fact, we are in a  case where  the operations have no symmetry, and the relations
leave the variables in the same order. So the operad is \emph {regular}, that is it is
determined by a non-$\Sigma$-operad: $\PP_n = {\mathbb K}[X_{n-1}]$. 

Reversing the roles of $X$ and $Z$, that is taking the elements of $X$ as relations, gives
rise to a new (non-$\Sigma$-)operad $\mathcal Q$ such that ${\mathcal Q}_n={\mathbb K}[Z_{n-1}]$.

\begin{lemma} \label{thedual}
The Koszul dual operad of $\PP$ is $\mathcal Q$, that is $\PP^! = \mathcal Q .$
\end{lemma}

\noindent Proof. Recall from \cite{G-K}, (see \cite{L1} for a short survey and \cite{Fresse} for details) that the dual operad $\PP^!$ of the
non-$\Sigma$-operad $\PP$ is constructed as follows. The generating operations are the
same. The space of relations is made of the elements $\sum \alpha_{ij} (x\circ_i y)\circ_j z
+ \sum \beta_{ij} x\circ_i ( y\circ_j z )$ (for some scalars $\alpha_{ij}$ and $\beta_{ij}$) which are orthogonal to the relations of $\PP$ for
the inner product $\langle -,-\rangle$ defined on the linear generators by
\begin{eqnarray}
\langle (x\circ_i y)\circ_j z\ ,\ (x\circ_i y)\circ_j z\rangle = 1, \\
\langle x\circ_i ( y\circ_j z )\ ,\ x\circ_i ( y\circ_j z )\rangle= -1 \\
\langle -,-\rangle= 0\qquad  {\rm otherwise.}
\end{eqnarray}

One immediately checks that the vector space generated by $X$ is orthogonal to the vector
space generated by $Z$, and therefore the Koszul dual of $\PP$ is $\mathcal Q$. \hfill $\square$

\begin{thm}\label{duality}
The operads $\PP$ and $\mathcal Q$ are Koszul operads.
\end{thm}
\noindent Proof. The Koszul duality of $\PP$ is equivalent to the acyclicity of the Koszul complex of
$\PP$, which is $(\PP^{!*}(\PP(V)), \delta)$. Since $\PP$ is regular (i.e.~comes from a
non-$\Sigma$-operad), it is sufficient to check the acyclicity for $V={\mathbb K}$.  Since
$\PP^!={\mathcal Q}$ the chains of the Koszul complex of $\PP$ are the same as the chains
of the Koszul complex of $X$ constructed in the first section. A careful checking of the
construction of $\delta$ shows that $\delta = d$.

So we can apply Corollary~\ref{acyclicity} and the proof is completed.\hfill $\square$ 

\subsection{Remarks}
 The Poincar\'e series of an operad is defined as 
$$f^{\PP}(t) := \sum _{n\geq 1} (-1)^n \frac {\dim \PP(n)}{n!}t^n = \sum _{n\geq 1} (-1)^n
\dim \PP_n\ t^n\ .$$
Hence, for the operad $\PP$ defined by $X$, one has $f^{\PP}(t)= f(X,t)$ and the functional equation of Theorem \ref{THM} is the functional equation 
$f^{\PP^!}(f^{\PP}(t))= t$ proved in \cite{G-K} for Koszul operads.

In this paper we exploit only the Poincar\'e series property of Koszul operads. There are many other applications like constructing homotopy algebras (cf.~\cite{G-K}) and computing the homology of the associated partition complex (cf.~\cite{V2}).

\section{Generalization}
There is no reason to restrict oneself to binary trees, that is to binary operads. One can
start with planar rooted trees. In this framework we choose a set of index for each integer $k\geq 2$. Hence the functional equation is now in two variables, see~\cite{V2} for the operadic interpretation.  In this section we give some examples of a particular case: the vertices of the trees have valence $k\geq 2$ for a fixed $k$.

\subsection{$k$-ary planar trees}
Let $Y_n^{(k)}$ be the set of planar rooted trees with $n$ vertices, each vertex being of valence $k$. The number of
leaves of such a tree is $(k-1)n+1$. The case $k=2$ is the one treated in the first part. Let $I$ be a set of indices and
let $Y_n^{(k)}\times I^n$ be the set of labelled trees. Choose a subset $X$ of $Y_2^{(k)}\times I^2$ and let $Z$ be its
complement. As before we define $X_n\subset Y_n^{(k)}\times I^n$ to be the subset made of labelled trees whose
local patterns belong to $X$.

In order to state the Theorem we need to introduce the following series. Let ${\underline a} = (a_0, \ldots,
a_n,\ldots )$ be a sequence of numbers (we will always have $a_0=1$). Define the (lacunary) series $f^{(k)}$ and
$g^{(k)}$ as follows:
$$
f^{(k)}({\underline a}, t) := \sum _{n\geq 0}(-1)^{n+1} a_n t^{(k-1)n+1} = -t + a_1t^k-a_2t^{2k-1} +\cdots 
$$
$$
g^{(k)}({\underline a}, t) := - \sum _{n\geq 0}(-1)^{(k+1)n} a_n t^{(k-1)n+1} = -t +
(-1)^ka_1t^k- a_2t^{2k-1} + \cdots
$$
Observe that when $k$ is even $f^{(k)}= g^{(k)}$ and when $k$ is odd all the signs in $g^{(k)}$ are $-$ .
When $k=2$, one has $f^{(2)}= g^{(2)}= f$ as defined in section 2. The series $f^{(k)}(X,t)$ and $g^{(k)}(X,t)$ are obtained
by taking $a_n= \#X_n$.

\begin{thm}\label{THMdeux}  Let $X$ be a subset of  $Y_2^{(k)}\times I^2$ and let $Z$ be its
complement, i.e.~$X \sqcup Z= Y_2^{(k)}\times I^2$. Then the following functional equation holds:

$$g^{(k)}(Z,f^{(k)}(X,t)) = t .$$
\end{thm}
The proof is along the same line as the proof of Theorem~\ref{THM} and we let the diligent reader to verify it.

There is an operadic interpretation of this result, which involves the notion of $k$-ary algebras. The
relevant generalization of Koszul duality theory for quadratic algebras (not just binary) can be found in \cite{Fresse}.

It would be interesting to study the analogous question with operads replaced by props as in \cite{V1}.

\subsection{Examples}
The integer sequences involved in this case are of the form 
\begin{displaymath}
(1, \underbrace{0, \ldots, 0}_{k-2}, a_1, \underbrace{0, \ldots, 0}_{k-2}, a_2, \underbrace{0, \ldots, 0}_{k-2}, a_3,
0, \ldots, )
\end{displaymath}
with $a_1 = \# I, a_2= \# X$, so  $0\leq a_2\leq k(a_1)^2$. We denote such a lacunary sequence by 
$(1; a_1; a_2; \cdots ; a_n; \cdots )_k$ .
\\

\noindent (a) $(1;  1; k;\frac{k(3k-1)}{2}; \frac{k(8k^2-6k+1)}{3}; \cdots )_k$ versus $(1; 1; 0; \cdots ; 0; \cdots )_k$ .

Let $c_n^{(k)}$ be the number of $k$-ary trees with $n$ vertices. Taking $I=\{1\}$, $X=Y_2^{(k)}$ and $Z=\emptyset$ we get $g^{(k)}(\emptyset,t)=-t+(-1)^kt^k$ and $f^{(k)}(Y_2^{(k)},t)=\sum_{n\geq 0}(-1)^{n+1}c_n^{(k)}t^{(k-1)n+1}$. So this last series, denote it $y$, satisfies the functional equation $-y+ (-1)^ky^k = t$.
\\

\noindent (b)  $(1;  1; 2; 5; \cdots ; c_n ; \cdots )_3$ versus  $(1;  1; 1; \cdots ; 1; \cdots )_3$ .

Take $I={1}$, $X$ has two elements and $Z$ has one element. The set $X_n$ has $c_{n+1}$ elements and $Z_n$ has one element. This case is related to the notion of totally associative ternary algebras and partially  associative ternary algebras studied in \cite {Gnedbaye}.
\\

\noindent (c) $(1;  1; 2; 5; \cdots ; c_n ; \cdots )_4$ versus itself.

Let $k=4$,  $I=\{1\}$. The set $X$ is made of two elements of $Y_2^{(4)}$ and $Z$ is made of the other two. It is clear that the sets $X_n$ and $Z_n$ are in bijection with the planar binary trees of degree $n$. As a consequence of Theorem~\ref{THMdeux} the series $h(t) = \sum_{n\geq 0}(-1)^{n+1}c_nt^{3n+1}$ satisfies $$h(h(t)) = t.$$
 Of course this result can also be proved by direct computation from the expression 
$c(t):= \sum_{n\geq 0}c_nt^{n}=\frac{1-\sqrt{1-4t}}{2t}$. 


\end{document}